\documentclass{article}
\usepackage{amssymb}
\usepackage{amsfonts}
\usepackage{amsmath}
\usepackage{geometry}
\usepackage{hyperref}

\newtheorem{theorem}{Theorem}

\newtheorem{proposition}[theorem]{Proposition}
\newtheorem{remark}[theorem]{Remark}

\newenvironment{proof}[1][Proof]{\noindent\textbf{#1.} }{\ \rule{0.5em}{0.5em}}
\input{tcilatex}
\begin{document}

\title{A generalized Fernique theorem and applications}
\author{Peter Friz\thanks{%
TU Berlin und WIAS. Email: friz@math.tu-berlin.de} \and Harald Oberhauser%
\thanks{%
Corresponding Author. Department of Pure Mathematics and Mathematical
Statistics, University of Cambridge. Email h.oberhauser@statslab.cam.ac.uk}\ 
\thanks{%
Supported by EPSCR Grant EP/P502365/1 and a DOC-fellowship of the Austrian
Academy of Sciences.}}
\maketitle

\begin{abstract}
We prove a generalisation of Fernique's theorem which applies to a class of
(measurable) functionals on abstract Wiener spaces by using the
isoperimetric inequality. Our motivation comes from rough path theory where
one deals with iterated integrals of Gaussian processes (which are
generically not Gaussian). Gaussian integrability with explicitly given
constants for variation and H\"{o}lder norms of the (fractional) Brownian
rough path, Gaussian rough paths and the Banach space valued Wiener process
enhanced with its L\'{e}vy area [Ledoux, Lyons, Quian. "L\'{e}vy area of
Wiener processes in Banach spaces". Ann. Probab., 30(2):546--578, 2002] then
all follow from applying our main theorem.
\end{abstract}

\section{A Generalized Fernique Theorem}

Let $\left( E,\left\vert \cdot \right\vert \right) $ be a real, separable
Banach space equipped with Borel $\sigma $-algebra $\mathfrak{B}$ and a
centered Gaussian measure $\mu $. A famous result by X.\ Fernique states
that $\left\vert \cdot \right\vert _{\ast }\mu $ has a Gauss tail; more
precisely, 
\begin{equation*}
\int \exp \left( \eta \left\vert x\right\vert ^{2}\right) \mathrm{d}\mu
\left( x\right) <\infty \text{ if }\eta <\frac{1}{2\sigma ^{2}},
\end{equation*}%
where 
\begin{equation}
\sigma :=\sup_{\xi \in E^{\ast },\left\vert \xi \right\vert _{E^{\ast
}}=1}\left( \,\int \left\langle \xi ,x\right\rangle ^{2}\mathrm{d}\mu \left(
x\right) \right) ^{1/2}<\infty ,  \label{DefSigmaAbstractWienerSpace}
\end{equation}%
and this condition on $\eta $ is sharp (see \cite[Thm 4.1]%
{Ledoux96:IsoperimetryandGaussianAnalysis} for instance). We recall the
notion of a \textit{reproducing kernel Hilbert space} $H$, continuously
embedded in $E$, $\left\vert h\right\vert \leq \sigma \left\vert
h\right\vert _{H}$ for all $h\in H$, so that $\left( E,H,\mu \right) $ is an
abstract Wiener space in the sense of L.\ Gross. We can then cite Borell's
inequality, e.g.\ \cite[Theorem 4.3]%
{Ledoux96:IsoperimetryandGaussianAnalysis}.

\begin{theorem}
\label{ThBorellTIS}Let $\left( E,H,\mu \right) $ be an abstract Wiener space
and $A\subset E$ a Borel set with $\mu \left( A\right) >0$. Take $a\in
(-\infty ,\infty ]$ such that%
\begin{equation*}
\mu \left( A\right) =\int_{-\infty }^{a}\frac{1}{\sqrt{2\pi }}e^{-x^{2}/2}%
\mathrm{d}x=:\Phi \left( a\right) .
\end{equation*}%
Then, if $\mathcal{K}$ denotes the unit ball in $H$ and $\mu _{\ast }$
stands for the inner measure\footnote{%
Measurability of the so-called Minkowski sum $A+r\mathcal{K}$ is a delicate
topic. Use of the inner measure bypasses this issue and is not restrictive
in applications.} associated to $\mu $,%
\begin{equation}
\mu _{\ast }\left( A+r\mathcal{K}\right) =\mu _{\ast }\left\{ x+rh:x\in
A,\,h\in \mathcal{K}\right\} \geq \Phi \left( a+r\right) .
\label{EqBorellInequality}
\end{equation}
\end{theorem}

The reader should observe that the following theorem reduces to the usual
Fernique result when applied to the Banach norm on $E$.

\begin{theorem}[Generalized Fernique]
\label{ThGeneralizedFernique}Let $\left( E,H,\mu \right) $ be an abstract
Wiener space. Assume $f:E\rightarrow \mathbb{R\cup }\left\{ -\infty ,\infty
\right\} $ is a measurable map and $N\subset E$ a null-set and $c$ some
positive constant such that for all $x\notin N$ 
\begin{equation}
\left\vert f\left( x\right) \right\vert <\infty  \label{EqfFinite}
\end{equation}%
\begin{equation}
\forall h\in H\text{: }\left\vert f\left( x\right) \right\vert \leq c\left\{
\left\vert \left( f\left( x-h\right) \right) \right\vert +\sigma \left\vert
h\right\vert _{H}\right\} .  \label{Control_H_translate}
\end{equation}%
Then, with the definition of $\sigma $ given in (\ref%
{DefSigmaAbstractWienerSpace}),%
\begin{equation*}
\int \exp \left( \eta \left\vert f\left( x\right) \right\vert ^{2}\right) 
\mathrm{d}\mu \left( x\right) <\infty \text{ \ if }\eta <\frac{1}{%
2c^{2}\sigma ^{2}}.
\end{equation*}
\end{theorem}

\begin{proof}
We have for all $x\notin N$\ and all $h\in r\mathcal{K}$, where $\mathcal{K}$
denotes the unit ball of $H$ and $r>0,$%
\begin{eqnarray*}
\left\{ x:\left\vert f\left( x\right) \right\vert \leq M\right\} &\supset
&\left\{ x:c\left( \left\vert f\left( x-h\right) \right\vert +\sigma
\left\vert h\right\vert _{H}\right) \leq M\right\} \\
&\supset &\left\{ x:c\left( \left\vert f\left( x-h\right) \right\vert
+\sigma r\right) \leq M\right\} \\
&=&\left\{ x+h:\left\vert f\left( x\right) \right\vert \leq M/c-\sigma
r\right\} \text{.}
\end{eqnarray*}%
Since $h\in r\mathcal{K}$ was arbitrary,%
\begin{eqnarray*}
\left\{ x:\left\vert f\left( x\right) \right\vert \leq M\right\} &\supset
&\cup _{h\in r\mathcal{K}}\left\{ x+h:\left\vert f\left( x\right)
\right\vert \leq M/c-\sigma r\right\} \\
&=&\left\{ x:\left\vert f\left( x\right) \right\vert \leq M/c-\sigma
r\right\} +r\mathcal{K}
\end{eqnarray*}%
and we see that%
\begin{eqnarray*}
\mu \left[ \left\vert f\left( x\right) \right\vert \leq M\right] &=&\mu
_{\ast }\left[ \left\vert f\left( x\right) \right\vert \leq M\right] \\
&\geq &\mu _{\ast }\left( \left\{ x:\left\vert f\left( x\right) \right\vert
\leq M/c-\sigma r\right\} +r\mathcal{K}\right) .
\end{eqnarray*}%
We can take $M=\left( 1+\varepsilon \right) c\sigma r$ and obtain%
\begin{equation*}
\mu \left[ \left\vert f\left( x\right) \right\vert \leq \left( 1+\varepsilon
\right) c\sigma r\right] \geq \mu _{\ast }\left( \left\{ x:\left\vert
f\left( x\right) \right\vert \leq \varepsilon \sigma r\right\} +r\mathcal{K}%
\right) .
\end{equation*}%
Keeping $\varepsilon $ fixed, take $r\geq r_{0}$ where $r_{0}$ is chosen
large enough such that $\mu \left[ \left\{ x:\left\vert f\left( x\right)
\right\vert \leq \varepsilon \sigma r_{0}\right\} \right] >0.$ Letting $\Phi 
$ denote the distribution function of a standard Gaussian, it follows from
Borell's inequality\ that%
\begin{equation*}
\mu \left[ \left\vert f\left( x\right) \right\vert \leq \left( 1+\varepsilon
\right) c\sigma r\right] \geq \Phi \left( a+r\right)
\end{equation*}%
for some $a>-\infty $. Equivalently, 
\begin{equation*}
\mu \left[ \left\vert f\left( x\right) \right\vert \geq x\right] \leq \bar{%
\Phi}\left( a+\frac{x}{\left( 1+\varepsilon \right) c\sigma }\right)
\end{equation*}%
with $\bar{\Phi}\equiv 1-\Phi $ and using $\bar{\Phi}\left( z\right)
\lesssim \exp \left( -z^{2}/2\right) $ this we see that this implies 
\begin{equation*}
\int \exp \left( \eta \left\vert f\left( x\right) \right\vert ^{2}\right) 
\mathrm{d}\mu \left( x\right) <\infty
\end{equation*}%
provided $\eta <\frac{1}{2}\left( \frac{1}{\left( 1+\varepsilon \right)
c\sigma }\right) ^{2}.$ Sending $\varepsilon \rightarrow 0$ finishes the
proof.
\end{proof}

\section{Applications\label{SecAppllicationsI}}

The examples below apply theorem \ref{ThGeneralizedFernique} with $f\left(
.\right) =\left\vert \left\vert \left\vert \cdot \right\vert \right\vert
\right\vert _{p\text{-var;}\left[ 0,T\right] }\circ \phi $ or $f\left(
.\right) =\left\vert \left\vert \left\vert \cdot \right\vert \right\vert
\right\vert _{1/p\text{-H\"{o}l;}\left[ 0,T\right] }\circ \phi $ where%
\begin{equation*}
\phi :C_{o}\left( \left[ 0,T\right] ,B\right) \rightarrow C_{o}\left( \left[
0,T\right] ,\left\{ 1\right\} \oplus B\mathbb{\oplus }B^{\otimes 2}\right)
\end{equation*}%
is constructed (typically through an almost-sure convergence result) as
measurable map, in general not continuous, and, for $\mathbf{x=}\left( 1,%
\mathbf{x}^{1}\mathbf{,x}^{2}\right) \,\mathbf{\in \,}C_{o}\left( \left[ 0,T%
\right] ,\left\{ 1\right\} \oplus B\mathbb{\oplus }B^{\otimes 2}\right) $,%
\begin{eqnarray*}
|||\mathbf{x|||}_{p\text{-var;}\left[ 0,T\right] } &=&\max_{i=1,2}\left(
\sup_{D=\left( t_{i}\right) }\sum_{i:t_{i}\in D}\left\vert \mathbf{x}%
_{t_{i},t_{i+1}}^{i}\right\vert _{B^{\otimes i}}^{p/i}\right) ^{1/p}\in %
\left[ 0,\infty \right] , \\
|||\mathbf{x|||}_{1/p\text{-H\"{o}l;}\left[ 0,T\right] }
&=&\max_{i=1,2}\left( \sup_{s,t\in \left[ 0,T\right] ,s\neq t}\frac{%
\left\vert \mathbf{x}_{s,t}^{i}\right\vert _{B^{\otimes i}}^{1/i}}{%
\left\vert t-s\right\vert ^{1/p}}\right) \in \left[ 0,\infty \right] .
\end{eqnarray*}%
This setting is standard in rough path theory (we refer to \cite%
{lyons-qian-02},\cite{lyons-caruana-levy-07}). The first example to have in
mind is $B\mathbb{=R}^{d}$ and $E=C_{o}\left( \left[ 0,T\right] ,\mathbb{R}%
^{d}\right) $ equipped with Wiener measure. We now discuss this in detail,
followed by more complex examples.

\subsection{Brownian motion on $\mathbb{R}^{d}$}

Take the usual Wiener space $\left( E,H,\mu \right) $, i.e.\ $E$ the space
of continuous paths on $\left[ 0,T\right] $ in $\mathbb{R}^{d}$ vanishing at 
$0$, $H$ the standard Cameron-Martin space and $\mu $ the Wiener measure.
The coordinate process $B_{t}\left( x\right) =x\left( t\right) $ for $x\in E$
is a standard Brownian motion. In this case, 
\begin{equation*}
\phi :C_{o}\left( \left[ 0,T\right] ,\mathbb{R}^{d}\right) \rightarrow
C_{o}\left( \left[ 0,T\right] ,G^{2}\left( \mathbb{R}^{d}\right) \right) ,
\end{equation*}%
(we recall that $G^{2}\left( \mathbb{R}^{d}\right) \subset \left\{ 1\right\}
\oplus \mathbb{R}^{d}\mathbb{\oplus }\left( \mathbb{R}^{d}\right) ^{\otimes
2}$ is the natural state space for geometric rough paths) is given by 
\begin{equation*}
\phi \left( x\right) =\left( 1,x,\int_{0}^{\cdot }x\otimes \circ dx\right)
\end{equation*}%
and almost surely well-defined (the final integral is a stochastic integral
in Stratonovich sense). We also write $\mathbf{B}=\left( 1,\mathbf{B}^{1},%
\mathbf{B}^{2}\right) $ for the coordinate process of the lift $\phi \left(
.\right) $ and call $\mathbf{B}$ enhanced Brownian motion. Using theorem \ref%
{ThGeneralizedFernique} it is easy to see that $\left\vert \left\vert
\left\vert \mathbf{B}\right\vert \right\vert \right\vert _{1/p\text{-H\"{o}l;%
}\left[ 0,T\right] }$ has a Gauss tail.

\begin{proposition}
\label{PropBM}If $p>2$%
\begin{equation*}
\mathbb{E}\left[ \exp \left( \frac{\eta }{T^{1-2/p}}\left\vert \left\vert
\left\vert \mathbf{B}\right\vert \right\vert \right\vert _{1/p\text{-H\"{o}l;%
}\left[ 0,T\right] }^{2}\right) \right] <\infty
\end{equation*}%
for all $\eta <\eta _{0}$. Moreover, one can take $\eta _{0}=\frac{1}{\left(
1+\sqrt{2}\right) }$.
\end{proposition}

\begin{proof}
We recall that $\left\vert \left\vert \left\vert \mathbf{B}\right\vert
\right\vert \right\vert _{p\text{-var;}\left[ 0,T\right] }<\infty $ a.s.\
(see \cite{lyons-qian-02}) and the stronger statement $\left\vert \left\vert
\left\vert \mathbf{B}\right\vert \right\vert \right\vert _{1/p\text{-H\"{o}l;%
}\left[ 0,T\right] }<\infty $ a.s.\ is found in \cite%
{Friz05:HoelderRoughPaths}. Since $\left\vert \left\vert \left\vert \mathbf{B%
}\right\vert \right\vert \right\vert _{1/p\text{-H\"{o}l;}\left[ 0,T\right]
}\sim T^{1/2-1/p}\left\vert \left\vert \left\vert \mathbf{B}\right\vert
\right\vert \right\vert _{1/p\text{-H\"{o}l;}\left[ 0,1\right] }$ we can
assume w.l.o.g.\ that $T=1$ and one checks (\ref{Control_H_translate}) by
simple Riemann--Stieltjes estimates. Note that by construction of $\mathbf{B}
$ (as an a.s.\ limit) and continuity properties of the integrals\footnote{%
The integral $\int_{s}^{t}h_{s,r}\otimes dx_{r}$ is defined as the
Riemann-Stieltjes integral $h_{s,t}\otimes x_{t}-\int_{s}^{t}dh_{r}\otimes
x_{r}$. The other integrals make immediate sense as Riemann-Stieltjes
integrals due to $\left\vert h\right\vert _{1-\text{var}}<\infty .$} it
follows that the set 
\begin{equation*}
\left\{ x\in E:\mathbf{B}\left( x+h\right) =\left( x+h\right) +\left( 
\mathbf{B}^{2}\left( x\right) +\int x\otimes dh+\int h\otimes dx+\int
h\otimes dh\text{ }\right) \text{ for all }h\in H\right\}
\end{equation*}%
has full measure (see \cite{CassFrizVictoir:WienerFunct}). Hence, there
exists a nullset $N$ s.t.\ for $x\notin N$ 
\begin{eqnarray*}
\left\vert \mathbf{B}^{1}\left( x+h\right) _{s,t}\right\vert &=&\left\vert
x_{s,t}+h_{s,t}\right\vert \leq \left\vert \left\vert \left\vert \mathbf{B}%
\left( x\right) \right\vert \right\vert \right\vert _{1/p\text{-H\"{o}l;}%
\left[ 0,1\right] }\left( t-s\right) ^{1/p}+\left\vert h_{s,t}\right\vert \\
\sqrt{\left\vert \mathbf{B}^{2}\left( x+h\right) _{s,t}\right\vert } &=&%
\sqrt{\left\vert \mathbf{B}_{s,t}^{2}+\int_{s}^{t}x_{s,r}\otimes
dh_{r}+\int_{s}^{t}h_{s,r}\otimes dx_{r}+\int_{s}^{t}h_{s,r}\otimes
dh_{r}\right\vert }
\end{eqnarray*}%
for all $h\in H$. Now,%
\begin{equation*}
\sqrt{\left\vert \int_{s}^{t}x_{s,r}\otimes
dh_{r}+\int_{s}^{t}h_{s,r}\otimes dx_{r}\right\vert }\leq \sqrt{\left\vert
h\right\vert _{1-\text{var;}\left[ s,t\right] }\left\vert x_{s,t}\right\vert 
}\leq \frac{1}{\sqrt{2}}\left( \left\vert h\right\vert _{1\text{-var;}\left[
s,t\right] }+\left\vert x_{s,t}\right\vert \right)
\end{equation*}%
and by Cauchy--Schwarz $\left\vert h\right\vert _{1-var;\left[ s,t\right]
}\leq \left\vert t-s\right\vert ^{1/2}\left\vert h\right\vert _{H}$ which
implies $\sqrt{\left\vert \int_{s}^{t}h_{s,r}\otimes dh_{r}\right\vert }\leq
\left\vert t-s\right\vert ^{1/2}\left\vert h\right\vert _{H}$. Combining
these estimates leads to%
\begin{eqnarray*}
\frac{\sqrt{\left\vert \mathbf{B}^{2}\left( x+h\right) _{s,t}\right\vert }}{%
\left( t-s\right) ^{1/p}} &\leq &\left( 1+1/\sqrt{2}\right) \left(
\left\vert \left\vert \left\vert \mathbf{B}\right\vert \right\vert
\right\vert _{1/p\text{-H\"{o}l;}\left[ 0,1\right] }+\left( t-s\right)
^{1/2-1/p}\left\vert h\right\vert _{H}\right) \\
&\leq &\left( 1+1/\sqrt{2}\right) \left( \left\vert \left\vert \left\vert 
\mathbf{B}\right\vert \right\vert \right\vert _{1/p\text{-H\"{o}l;}\left[ 0,1%
\right] }+\left\vert h\right\vert _{H}\right)
\end{eqnarray*}%
and $\left( \ref{Control_H_translate}\right) $ holds for $\left\vert f\left(
.\right) \right\vert =\left\vert \left\vert \left\vert \mathbf{B}\left(
.\right) \right\vert \right\vert \right\vert _{1/p\text{-H\"{o}l};\left[ 0,1%
\right] }$ with $c=\left( 1+1/\sqrt{2}\right) $ and $\sigma =\sqrt{E\left[
B_{1}^{2}\right] }=1$.
\end{proof}

\begin{remark}
This implies an (exponential) integrability of L\'{e}vy area which cannot be
obtained by integrability properties of the second Wiener-It\^{o} chaos due
to the non-linearity of area increments, i.e. $A_{s,t}\neq A_{0,t}-A_{0,s}$.
\end{remark}

\begin{remark}
It is well known that there exists $C=C\left( d\right) $ such that%
\begin{equation*}
C^{-1}\left\Vert \mathbf{x}\right\Vert \leq |||\mathbf{x|||}\leq C\left\Vert 
\mathbf{x}\right\Vert
\end{equation*}%
where $\mathbf{x=}\left( 1,\mathbf{x}^{1},\mathbf{x}^{2}\right) $, $|||%
\mathbf{x|||}=\max \left( \left\vert \mathbf{x}^{1}\right\vert ,\left\vert 
\mathbf{x}^{2}\right\vert \right) $ and $\left\Vert .\right\Vert $ denotes
the Carnot-Carath\'{e}odory norm on the step-$2$ free nilpotent group with $%
d $ generators $G^{2}\left( \mathbb{R}^{d}\right) $. As a consequence%
\begin{equation*}
C^{-1}\left\Vert \mathbf{x}\right\Vert _{p\text{-var;}\left[ 0,T\right]
}\leq |||\mathbf{x|||}_{p\text{-var;}\left[ 0,T\right] }\leq C\left\Vert 
\mathbf{x}\right\Vert _{p\text{-var;}\left[ 0,T\right] }
\end{equation*}%
where $\left\Vert \mathbf{x}\right\Vert _{p\text{-var;}\left[ 0,T\right]
}=\left( \sup_{D=\left( t_{i}\right) }\sum_{i:t_{i}\in D}\left\Vert \mathbf{x%
}_{t_{i},t_{i+1}}\right\Vert ^{p}\right) ^{1/p}$ and a similar estimate
holds for $\left\Vert \mathbf{x}\right\Vert _{1/p\text{-H\"{o}l;}\left[ 0,T%
\right] }=\sup_{s,t\in \left[ 0,T\right] ,s\neq t}\frac{\left\Vert \mathbf{x}%
_{s,t}\right\Vert }{\left\vert t-s\right\vert ^{1/p}}.$ Hence, the Gauss
tail of $|||\mathbf{B|||}_{1/p\text{-H\"{o}l;}\left[ 0,T\right] }$ is
consistent with the known Gauss tail of $\left\Vert \mathbf{B}\right\Vert
_{1/p\text{-H\"{o}l;}\left[ 0,T\right] }$ as obtained in \cite%
{FrizVictoir04:NoteonGeomRoughPaths} by a precise tracking of constants in
the Garsia-Rodemich-Rumsey estimate.
\end{remark}

\begin{remark}
Optimal variation of Brownian motion is not measured in $p-$variation norm
but in $\psi -$variation with $\psi \left( x\right) =x^{2}/\max \left( \log
\log 1/x,1\right) $. More precisely, in \cite{Taylor72} it is established
that%
\begin{equation*}
\sup_{D=\left( t_{i}\right) }\sum_{i:t_{i}\in D}\psi \left( \left\vert
B_{t_{i+1}}-B_{t_{i}}\right\vert \right) <\infty \text{ a.s.}
\end{equation*}%
This gives rise to a $\psi -$variation norm $\left\vert B\right\vert _{\psi 
\text{-var}}$ and likewise, one can show (see \cite%
{Davie06:DiscreteApproxRoughPaths}) that the rough path $\mathbf{B}$ has
finite $\psi $-variation $\left\Vert \mathbf{B}\right\Vert _{\psi -var}$,
which is the optimal variational regularity enjoyed by $\mathbf{B}$. This is
important as it allows to solve rough differential equations driven by $%
\mathbf{B}$ under minimal regularity assumptions on vector fields. It is
then interesting to know that the generalized Fernique estimate can be used
to see that the random variable $\left\Vert \mathbf{B}\right\Vert _{\psi
-var}$ also has a Gauss-tail; a fact which would be difficult to obtain from
tracking constants in the Garsia-Rodemich-Rumsey estimate. (If one aims for L%
\'{e}vy-modulus, the optimal modulus regularity enjoyed by $\mathbf{B}$ it
is in contrary possible, to obtain Gauss tail by tracking GRR constants, 
\cite{FrizVictoir05:ApproximationsEBM}).
\end{remark}

\subsection{Gaussian processes on $\mathbb{R}^{d}$}

Let us generalize from Brownian to a $d$-dimensional, continuous, centered
Gaussian process $X=\left( X^{1},\dots ,X^{d}\right) $ with independent
components, assuming the covariance of $X$ to be of finite $\rho $-variation
in $2$D-sense, $\left\vert R\right\vert _{\rho \text{-var;}\left[ 0,T\right]
}<\infty $ for some $\rho \in \left[ 1,2\right) $ (as introduced in \cite%
{FrizVictoir07:GaussRP1}; this setting covers for instance fractional
Brownian motion where $\rho =\frac{1}{2H}$ for a Hurst parameter $H\in
\left( \frac{1}{4},\frac{1}{2}\right] $, Ornstein-Uhlenbeck process, etc.).
The setup is as in the Brownian case, we just replace Wiener measure by a
more general Gaussian measure $\mu $ and an appropriate Cameron-Martin space 
$H$ and make the assumption of complementary Young regularity, i.e.\ $%
\exists q:1/q+1/\left( 2\rho \right) >1$ s.t.\ $H\hookrightarrow C^{q\text{%
-var}}$ (this assumption is always satisfied when $\rho \in \left[
1,3/2\right) $ and hence includes fBM with $H\in \left( \frac{1}{3},\frac{1}{%
2}\right] $).

\begin{proposition}
Let $X$ be a centered, continuous Gaussian process in $\mathbb{R}^{d}$ on $%
\left[ 0,T\right] $ with independent components and with covariance of
finite $2$D $\rho $-variation $\left\vert R\right\vert _{\rho \text{-var;}%
\left[ 0,T\right] }$ for $\rho <2$. Assume furthermore complementary Young
regularity. Then there exists a lift to a Gaussian rough path $\mathbf{X}%
\left( .\right) \in C_{0}\left( \left[ 0,T\right] ,G^{2}\left( \mathbb{R}%
^{d}\right) \right) $ of finite homogeneous $p$-variation, $p>2\rho $ and $%
|||\mathbf{X|||}_{p\text{-var;}\left[ 0,T\right] }$ has a Gauss tail.\ More
precisely, 
\begin{equation*}
\mathbb{E}\left[ \exp \left( \eta \frac{1}{\left\vert R\right\vert _{\rho 
\text{-var};\left[ 0,T\right] ^{2}}}|||\mathbf{X|||}_{p\text{-var;}\left[ 0,T%
\right] }^{2}\right) \right] <\infty
\end{equation*}%
for every $\eta <\eta _{0}$. Moreover, one can take $\eta _{0}=\left( \sqrt{2%
}3^{\left( 1/2-1/p\right) }\left( \sqrt{c_{\rho ,\rho }}^{p}+\sqrt{c_{\rho
,p}}^{p}/\sqrt{2}\right) ^{1/p}\right) ^{-2}$ where $c_{u,v}=2.4^{1/u+1/v}%
\zeta \left( \frac{1}{u}+\frac{1}{v}\right) $ (with the Riemann-Zeta
function $\zeta \left( s\right) =\sum_{n=1}^{\infty }n^{-s}$ ).
\end{proposition}

\begin{proof}
The lift constructed in \cite{FrizVictoir07:GaussRP1} gives a rough path of
finite $p$-variation and we verify the translation estimate $\left( \ref%
{Control_H_translate}\right) $. By construction of the rough path lift (as
an a.s. limit) and the assumption of complementary Young regularity one has
that 
\begin{equation*}
\left\{ x:\mathbf{X}\left( x+h\right) =\left( x+h\right) +\left( \mathbf{X}%
^{2}\left( x\right) +\int h\otimes dx+\int x\otimes dh+\int h\otimes
dh\right) \text{ for all }h\in H\right\}
\end{equation*}%
has full measure. The estimate $\left( a+b\right) ^{p}\leq 2^{p-1}\left(
a^{p}+b^{p}\right) $ together with $\left\vert h\right\vert _{\rho \text{%
-var;}\left[ s,t\right] }\leq \left\vert h\right\vert _{H}\sqrt{\left\vert
R\right\vert _{\rho \text{-var;}\left[ s,t\right] ^{2}}}$ (see \cite%
{FrizVictoir07:GaussRP1}) implies 
\begin{equation*}
\left\vert \mathbf{X}^{1}\left( x+h\right) _{s,t}\right\vert ^{p}\leq
2^{p-1}\left( \left\vert x_{s,t}\right\vert ^{p}+\left\vert h\right\vert
_{H}^{p}\left( \left\vert R\right\vert _{\rho \text{-var};\left[ s,t\right]
^{2}}^{\rho }\right) ^{p/\left( 2\rho \right) }\right)
\end{equation*}%
and since by assumption $p/\left( 2\rho \right) >1$ and $\left\vert
R\right\vert _{\rho \text{-var};\left[ 0,T\right] ^{2}}^{\rho }$ is a $2$%
D-control summing up yields 
\begin{equation*}
\sum_{i}\left\vert \mathbf{X}^{1}\left( x+h\right)
_{t_{i},t_{i+1}}\right\vert _{B}^{p}\leq 2^{\left( p-1\right) /p}\left(
\left\vert \left\vert \left\vert \mathbf{X}\left( x\right) \right\vert
\right\vert \right\vert _{p\text{-var;}\left[ 0,T\right] }^{p}+\left\vert
h\right\vert _{H}^{p}\left( \left\vert R\right\vert _{\rho \text{-var};\left[
0,T\right] ^{2}}^{\rho }\right) ^{p/\left( 2\rho \right) }\right) \text{.}
\end{equation*}%
Similar for the second level 
\begin{equation*}
\left\vert \mathbf{X}^{2}\left( x+h\right) _{s,t}\right\vert _{B^{\otimes
2}}^{p/2}\leq 3^{p/2-1}\left( \left\vert \mathbf{X}_{s,t}^{2}\left( x\right)
\right\vert _{B\otimes B}^{p/2}+\left\vert \int_{s}^{t}h_{s,r}\otimes
dx_{r}+\int_{s}^{t}x_{s,r}\otimes dh_{r}\right\vert _{B\otimes
B}^{p/2}+\left\vert \int_{s}^{t}h_{s,r}\otimes dh_{r}\right\vert _{B\otimes
B}^{p/2}\right) \text{.}
\end{equation*}%
By Young's inequality (and using an i.b.p.) there exists a constant $c_{\rho
,p}$ such that 
\begin{eqnarray*}
\left\vert \int_{s}^{t}x_{s,r}\otimes dh_{r}+\int_{s}^{t}h_{s,r}\otimes
dx_{r}\right\vert ^{p/2} &=&\left\vert \int_{s}^{t}dh_{r}\otimes
x_{s,t}\right\vert ^{p/2}\leq \left( c_{\rho ,p}\left\vert h\right\vert
_{\rho \text{-var};\left[ s,t\right] }\left\vert x\right\vert _{p\text{-var};%
\left[ s,t\right] }\right) ^{p/2} \\
&\leq &\frac{\sqrt{c_{\rho ,p}}^{p}}{\sqrt{2}}\left( \left\vert h\right\vert
_{H}^{p}\left( \left\vert R\right\vert _{\rho \text{-var};\left[ 0,T\right]
^{2}}^{\rho }\right) ^{p/\left( 2\rho \right) }+\left\vert x\right\vert _{p%
\text{-var};\left[ s,t\right] }^{p}\right) \text{.}
\end{eqnarray*}%
Further, 
\begin{equation*}
\left\vert \int_{s}^{t}h_{s,r}\otimes dh_{r}\right\vert _{B\otimes
B}^{p/2}\leq \sqrt{c_{\rho ,\rho }}^{p}\left\vert h\right\vert _{\rho \text{%
-var};\left[ s,t\right] }^{p}\leq \sqrt{c_{\rho ,\rho }}^{p}\left\vert
h\right\vert _{H}^{p}\left( \left\vert R\right\vert _{\rho \text{-var};\left[
0,T\right] ^{2}}^{\rho }\right) ^{p/2\rho }\text{.}
\end{equation*}%
Combining these estimates yields 
\begin{eqnarray*}
\sum_{{}}\left\vert \mathbf{X}^{2}\left( x+h\right)
_{t_{i},t_{i+1}}\right\vert _{B}^{p/2} &\leq &3^{\left( p/2-1\right) }\left( 
\sqrt{c_{\rho ,\rho }}^{p}+\sqrt{c_{\rho ,p}}^{p}/\sqrt{2}\right) \left(
\left\vert \left\vert \left\vert \mathbf{X}\left( x\right) \right\vert
\right\vert \right\vert _{p\text{-var;}\left[ 0,T\right] }^{p}+\left\vert
h\right\vert _{H}^{p}\left( \left\vert R\right\vert _{\rho \text{-var};\left[
0,T\right] ^{2}}^{\rho }\right) ^{p/\left( 2\rho \right) }\right) \\
&\leq &3^{\left( p/2-1\right) }\left( \sqrt{c_{\rho ,\rho }}^{p}+\sqrt{%
c_{\rho ,p}}^{p}/\sqrt{2}\right) \left( \frac{\sqrt{\left\vert R\right\vert
_{\rho \text{-var};\left[ 0,T\right] ^{2}}}}{\sigma }\right) ^{p}\left(
\left\vert \left\vert \left\vert \mathbf{X}\left( x\right) \right\vert
\right\vert \right\vert _{p\text{-var;}\left[ 0,T\right] }^{p}+\left\vert
h\right\vert _{H}^{p}\sigma ^{p}\right) .
\end{eqnarray*}%
Noting that $\sqrt{\left\vert R\right\vert _{\rho \text{-var};\left[ 0,T%
\right] ^{2}}}/\sigma \geq 1$ (which follows directly from the definition of 
$\left\vert R\right\vert _{\rho \text{-var}}$ and $\sigma $) the above
estimates imply that $\left( \ref{Control_H_translate}\right) $ holds with $%
c=3^{\left( 1/2-1/p\right) }\left( \sqrt{c_{\rho ,\rho }}^{p}+\sqrt{c_{\rho
,p}}^{p}/\sqrt{2}\right) ^{1/p}\sqrt{\left\vert R\right\vert _{\rho \text{%
-var};\left[ 0,T\right] ^{2}}}/\sigma $ where $c_{u,v}=2.4^{1/u+1/v}\zeta
\left( \frac{1}{u}+\frac{1}{v}\right) $ and $\sigma $ as in $\left( \ref%
{DefSigmaAbstractWienerSpace}\right) $.
\end{proof}

\begin{remark}
One can replace the variation norm by a stronger H\"{o}lder norm (if finite)
and follow the proof of Proposition \ref{PropBM}.\bigskip\ A sufficient
condition for H\"{o}lder regularity is that the covariance of $X$ satisfies
an appropriate H\"{o}lder condition (see \cite{FrizVictoir07:GaussRP1}).
\end{remark}

\begin{remark}
Explicit estimates for $\left\vert R\right\vert _{\rho \text{-var;}\left[ 0,T%
\right] }$ are sometimes known. For example for Brownian motion $\left\vert
R\right\vert _{1\text{-var;}\left[ 0,T\right] }=T$ \ (compare this to
Proposition \ref{PropBanachWiener} where one has the same scaling in $T$ but
a better constant $\eta _{0}$ due to the use of Riemann-Stieltjes instead of
Young integrals) and for fBM with Hurst parameter $H\in \left( 0,1/2\right] $%
, $\left\vert R\right\vert _{1/\left( 2H\right) \text{-var;}\left[ 0,T\right]
}\leq C\left( H\right) T$ $^{1/\left( 2H\right) }$ for a known constant $%
C\left( H\right) $ (see \cite{FrizVictoir07:GaussRP1})
\end{remark}

\subsection{Banach space valued Brownian motion}

Let $\left( B,G,\nu \right) $ be an abstract Wiener space. Then there exists
a stochastic process $X=\left( X_{t}\right) $ with continuous paths, taking
values in $B$ and such that $X_{0}=0$ a.s., the increments $X_{t}-X_{s}$ are
independent and $\xi \left( X_{t}-X_{s}\right) \sim N\left( 0,\left(
t-s\right) \left\vert \xi \right\vert _{G}^{2}\right) $ for $\xi \in B^{\ast
},$ $s<t$. This process is called Wiener process based on (the abstract
Wiener space) $\left( B,G,\nu \right) $. The associated L\'{e}vy area (and
lift to a rough path $\mathbf{X}$) was constructed in \cite%
{ledoux-lyons-qian-02} (the choice of a norm $\left\vert .\right\vert
_{B^{\otimes 2}}$ on $B^{\otimes 2}$ is a subtle issue and reflected in the
condition of exactness; the injective tensor norm is always exact; see \cite%
{ledoux-lyons-qian-02} for definition and further examples).

\begin{proposition}
\label{PropBanachWiener}Let $X$ be a Wiener process based on the abstract
Wiener space $\left( B,G,\nu \right) $ and let $\sigma
_{B}^{2}=\sup_{\left\vert \varphi \right\vert _{B^{\ast }\leq
1}}\int_{B}\varphi \left( x\right) ^{2}\mu \left( dx\right) $. If we
complete the algebraic tensor product $B\otimes B$ with a norm $\left\vert
.\right\vert _{B\otimes B}$ such that the pair $\left( \left\vert
.\right\vert _{B\otimes B},\nu \right) $ is exact then there exists a lift
to a geometric rough path $\mathbf{X=}\left( 1,\mathbf{X}^{1},\mathbf{X}%
^{2}\right) $ $\in C^{p\text{-var}}\left( \left[ 0,T\right] ,1\oplus B\oplus
B^{\otimes 2}\right) $, $p>2$ and $|||\mathbf{X|||}_{p\text{-var}}$ has a
Gauss tail. More precisely, 
\begin{equation*}
\mathbb{E}\left[ \exp \left( \eta \frac{1}{\sigma _{B}^{2}T}|||\mathbf{X|||}%
_{p\text{-var;}\left[ 0,T\right] }^{2}\right) \right] <\infty
\end{equation*}%
for every $\eta <\eta _{0}$. Moreover, one can take $\eta _{0}=\left( \sqrt{2%
}3^{\left( 1/2-1/p\right) }\left( 1+1/\sqrt{2}\right) ^{1/p}\right) ^{-2}$.
\end{proposition}

Let $E=C_{0}\left( \left[ 0,T\right] ,B\right) $ and set $\left\vert
x\right\vert _{E}=\sup_{t\in \left[ 0,T\right] }\left\vert x_{t}\right\vert $%
. Denote associated Borel $\sigma -$algebra $\mathcal{E}$ and define the
RKHS $H\subset E$ as $H=i^{\ast }\left( L^{2}\left( E,\mathcal{E},\mu
\right) \right) $ where $i^{\ast }:L^{2}\left( E,\mathcal{E},\mu ;\mathbb{R}%
\right) \rightarrow $ $E$ maps $\varphi $ to $\int_{E}x\varphi \left(
x\right) \mu \left( dx\right) $ ($\mu $ is the measure given by construction
of the Wiener process $\left( X_{t}\right) $). $\left( E,H,\mu \right) $ is
then an abstract Wiener space (following the setup of \cite%
{Ledoux96:IsoperimetryandGaussianAnalysis}). We first estimate the
regularity of elements in the Hilbert space $H$.

\begin{proposition}
\label{PropRegularityRKHS}For all $h\in H,$ $\left\vert h\right\vert _{1%
\text{-var;}\left[ s,t\right] }\leq \left\vert h\right\vert _{H}\sigma
_{B}\left\vert t-s\right\vert ^{1/2}$.
\end{proposition}

\begin{proof}
Let $\left( t_{j}\right) $ be a dissection of $\left[ s,t\right] $. By
definition $h=\int_{E}x\varphi \left( x\right) \mu \left( dx\right) $ for
some $\psi \in L^{2}\left( \mu \right) $ and we also write $h_{t}=\mathbb{E}%
\left( \psi \left( X\right) X_{t}\right) $. Then 
\begin{eqnarray*}
\sum_{j}\left\vert h_{t_{j},t_{j+1}}\right\vert _{B} &=&\sup_{\xi _{j}\in
B^{\ast }:\left\vert \xi _{j}\right\vert _{B^{\ast }}\leq 1}\sum_{j}\xi
_{j}h_{t_{j},t_{j+1}}=\sup \mathbb{E}\left[ \psi \left( x\right) \sum_{j}\xi
_{j}X_{t_{j},t_{j+1}}\right] \\
&\leq &\left\vert \psi \right\vert _{L^{2}\left( \mu \right) }\sup 
_{\substack{ \xi _{j}\in B^{\ast }:\left\vert \xi _{j}\right\vert _{B^{\ast
}}\leq 1  \\ \xi _{k}\in B^{\ast }:\left\vert \xi _{k}\right\vert _{B^{\ast
}}\leq 1}}\sqrt{\sum_{j,k}\mathbb{E}\left[ \xi _{j}X_{t_{j},t_{j+1}}\xi
_{k}X_{t_{k},t_{k+1}}\right] } \\
&=&\left\vert h\right\vert _{H}^{2}\sup_{\substack{ \xi _{j}\in B^{\ast
}:\left\vert \xi _{j}\right\vert _{B^{\ast }}\leq 1  \\ \xi _{k}\in B^{\ast
}:\left\vert \xi _{k}\right\vert _{B^{\ast }}\leq 1}}\sqrt{\sum_{j,k}\mathbb{%
E}\left[ \xi _{j}X_{t_{j},t_{j+1}}\xi _{k}X_{t_{k},t_{k+1}}\right] }
\end{eqnarray*}%
where we used Cauchy-Schwarz and the isometry $\left\langle i^{\ast }\left( 
\tilde{\psi}\right) ,i^{\ast }\left( \psi \right) \right\rangle
_{H}=\left\langle \tilde{\psi},\psi \right\rangle _{L^{2}\left( \mu \right)
} $. The characteristic property of the tensor product yields an isomorphism
between the space of (algebraic) linear\ functionals on $B\otimes B$, $%
L\left( B\otimes B,\mathbb{R}\right) $, and bilinear functionals on $B\times
B$, $BL\left( B\times B,\mathbb{R}\right) $. Since $\xi _{j}\left( .\right)
\xi _{k}\left( .\right) \in BL\left( B\times B,\mathbb{R}\right) $ 
\begin{equation*}
\xi _{j}\left( .\right) \xi _{k}\left( .\right) =\left( \xi _{j}\otimes \xi
_{k}\right) \left( .\otimes .\right) \text{ for some }\xi _{j}\otimes \xi
_{k}\in \text{ }L\left( B\otimes B,\mathbb{R}\right)
\end{equation*}%
and $\left\vert \xi _{j}\otimes \xi _{k}\right\vert \leq 1$ since $%
\left\vert \xi _{j}\right\vert \leq 1,\left\vert \xi _{j}\right\vert \leq 1$%
. So we have the estimate%
\begin{equation*}
\sum_{j}\left\vert h_{t_{j},t_{j+1}}\right\vert \leq \left\vert h\right\vert
_{H}^{2}\sqrt{\left( \sum_{j,k}\left\vert \mathbb{E}\left(
X_{t_{j},t_{j+1}}\otimes X_{t_{k},t_{k+1}}\right) \right\vert _{B\otimes
B}\right) \text{.}}
\end{equation*}%
Again for $\varphi \in L\left( B\otimes B,\mathbb{R}\right) $, $\varphi
\left( \mathbb{E}\left( X_{t_{j},t_{j+1}}\otimes X_{t_{k},t_{k+1}}\right)
\right) =\mathbb{E}\left( \tilde{\varphi}\left(
X_{t_{j},t_{j+1}},X_{t_{k},t_{k+1}}\right) \right) $ for some $\tilde{\varphi%
}\in BL\left( B\times B,\mathbb{R}\right) $ by linearity of $\varphi $.
Writing $\tilde{\varphi}\left( x,y\right) =\sum_{i}f_{i}\left( x\right)
g_{i}\left( y\right) $, $f_{i},g_{i}\in L\left( B,\mathbb{R}\right) $, and
using independence of increments in combination with Gaussianity and $%
\left\vert .\right\vert _{B\otimes B}=\sup \left\{ \varphi \left( .\right)
,\varphi \in L\left( B\otimes B,\mathbb{R}\right) ,\left\vert \varphi
\right\vert \leq 1\right\} $ gives $\sum_{j\neq k}\left\vert \mathbb{E}%
\left( X_{t_{j},t_{j+1}}\otimes X_{t_{k},t_{k+1}}\right) \right\vert
_{B\otimes B}=0$.

By the compatibility of the tensor norm (and again using the characteristic
property of tensor products) 
\begin{equation*}
\left\vert \mathbb{E}\left( X_{t_{j},t_{j+1}}\otimes
X_{t_{j},t_{j+1}}\right) \right\vert _{B\otimes B}\leq \mathbb{E}\left(
\left\vert X_{t_{j},t_{j+1}}\right\vert _{B}^{2}\right) .
\end{equation*}%
\ Now $X_{t_{j},t_{j+1}}$ has the same distribution as $\left(
t_{j+1}-t_{j}\right) X_{1}$ and since $X_{1}$ has distribution $\mu $, 
\begin{equation*}
\mathbb{E}\left( \left\vert X_{t_{j},t_{j+1}}\right\vert _{B}^{2}\right)
=\left( t_{j+1}-t_{j}\right) \int_{B}\left\vert x\right\vert _{B}^{2}\mu
\left( dx\right) .
\end{equation*}%
This gives $\sum_{j}\left\vert h_{t_{j},t_{j+1}}\right\vert \leq \left\vert
h\right\vert _{H}\sigma _{B}\sqrt{t-s}$.
\end{proof}

\begin{proof}[Proof of Proposition \protect\ref{PropBanachWiener}]
Since $|||\mathbf{X|||}_{p\text{-var;}\left[ 0,T\right] }<\infty $ a.s.\ is
shown in \cite{ledoux-lyons-qian-02} we only have to check $\left( \ref%
{Control_H_translate}\right) $. We recall the scaling $|||\mathbf{X|||}_{p%
\text{-var;}\left[ 0,T\right] }\sim \sqrt{T}|||\mathbf{X|||}_{p\text{-var;}%
\left[ 0,1\right] }$ and that by construction of the Lyons-Ledoux-Quian
rough path lift (as an a.s.\ limit) and continuity properties of Riemann
integrals the set 
\begin{equation*}
\left\{ x:\mathbf{X}\left( x+h\right) =\left( x+h\right) +\left( \mathbf{X}%
^{2}\left( x\right) +\int h\otimes dx+\int x\otimes dh+\int h\otimes
dh\right) \text{ for all }h\in H\right\}
\end{equation*}%
has full measure. The proof now works as in the finite dimensional case.
\end{proof}

\begin{remark}
The Gauss tail of $\left\vert \left\vert \left\vert \mathbf{X}\right\vert
\right\vert \right\vert _{p\text{-var;}\left[ 0,T\right] }$ was obtained in 
\cite{Inahama-JFA-06} by a careful tracking of the original estimates in 
\cite{ledoux-lyons-qian-02} though no explicit constant $\eta $ was given.
\end{remark}

\subsection{Integrability of higher iterated integrals}

The above applications imply Gaussian integrability of norms of twice
iterated integrals. To obtain integrability properties of $N$-times iterated
integrals, $N\geq 3,$ we recall a basic theorem of rough path theory (see 
\cite{lyons-caruana-levy-07}): a continuous $G^{\left[ p\right] }\left( 
\mathbb{R}^{d}\right) $-valued path $\mathbf{x}$ of finite $p$-variation
lifts for every $N\geq \left[ p\right] $ uniquely to a $G^{N}\left( \mathbb{R%
}^{d}\right) $-valued path, say $S_{N}\left( \mathbf{x}\right) $, of finite $%
p$-variation and there exists a constant $C\left( N,p\right) $ such that%
\begin{equation*}
\left\Vert S_{N}\left( \mathbf{x}\right) \right\Vert _{p\text{-var}}\leq
C\left( N,p\right) \left\Vert \mathbf{x}\right\Vert _{p\text{-var}}\text{.}
\end{equation*}%
Applied to, say, $d$-dimensional Brownian motion, this yields in combination
with proposition \ref{PropBM} that Brownian motion and all its iterated
Stratonovich integrals up to order $N,$ written as $S_{N}\left( \mathbf{B}%
\right) $ and viewed as a diffusion in the step-$N$ nilpotent group with $d$
generators have Gaussian integrability in the sense that $\left\Vert
S_{N}\left( \mathbf{B}\right) \right\Vert _{p\text{-var}}$ has a Gauss tail.

\textbf{Acknowledgement:} The authors are indebted to Nicolas Victoir, Terry
Lyons and James Norris for helpful conversations.

\bibliographystyle{plain}
\bibliography{rpath}

\end{document}